\documentclass[11pt, reqno]{amsart}
\def\squarebox#1{\hbox to #1{\hfill\vbox to #1{\vfill}}}

\newcommand{\R}{{\mathbb R}}

\newcommand{\N}{{\mathbb N}}

\theoremstyle{plain}

\newtheorem{thm}{Theorem}[section]

\newtheorem{lem}{Lemma}[section]
\newtheorem{rem}{Remark}
\newtheorem{prop}{Proposition}[section]

\renewcommand{\thelem}{\thesection.\arabic{lem}}

\renewcommand{\theequation}{\thesection.\arabic{equation}}
\renewcommand{\therem}{\thesection.\arabic{rem}}

\numberwithin{equation}{section}

\begin{document}

\title[Nonlinear wave equation]
{Polynomial bounds on the Sobolev norms of the solutions of the nonlinear wave equation with time dependent potential}

\author[V. Petkov]{Vesselin Petkov}
\address{Institut de Math\'ematiques de Bordeaux, 351,
Cours de la Lib\'eration, 33405  Talence, France}
\email{petkov@math.u-bordeaux.fr}
\author[N. Tzvetkov]{Nikolay Tzvetkov}
\address{D\'epartement de Math\'ematiques (AGM ), Universit\'e de Cergy-Pontoise, 2, av. Adolphe Chauvin, 95302 Cergy-Pontoise Cedex, France}
\email{nikolay.tzvetkov@u-cergy.fr}



\def\ts{\tilde{\sigma}}
\def\p{{\mathcal P}}
\def\h{{\mathcal H}}
\def\pa{\partial}
\def\om{\omega}
\def\fm{\frac{1}{1+m}}
\def\gy{\frac{1 + y}{\gamma}}
\def\nn{\Bigl( 1 - \frac{1}{n+1}\Bigr)}
\def\pa{\partial}
\def\ep{\epsilon}
\def\Bl{\Bigl (}
\def\Br{\Bigr )}
\def\r3{\R^3}
\def\l2{L^2(\r3)}
\def\Li6{L^{\infty}([s, s+ \tau], L^6(\r3))}
\def\hc{{\mathcal H}}

\keywords{Time periodic potential, Nonlinear wave equation,
 Growth of Sobolev norms}

\begin{abstract} We consider the Cauchy problem for the nonlinear wave equation $u_{tt} - \Delta_x u +q(t, x) u + u^3 = 0$ with smooth potential $q(t, x) \geq 0$ having compact support with respect to $x$. The linear equation without the nonlinear term $u^3$ and potential periodic in $t$ may have solutions with exponentially increasing as $ t \to \infty$ norm $H^1(\R^3_x)$. In \cite{PT} it was established that adding the nonlinear term $u^3$ the $H^1(\R^3_x)$ norm of the solution is polynomially bounded for every choice of $q$.  In this paper we show that $H^k(\R^3_x)$ norm of this global solution is also polynomially bounded. To prove this we apply a different argument based on the analysis of a sequence $\{Y_k(n\tau_k)\}_{n = 0}^{\infty}$ with suitably defined energy norm $Y_k(t)$ and $0 < \tau_k <1.$
\end{abstract}

 \maketitle


\section{Introduction}

Consider for $t \in \R,\: x \in \R^3$ the Cauchy problem 
\begin{equation} \label{eq:1.1}
\pa_t^2 u - \Delta_x u + q(t, x) u + u^3 = 0, \: u(0, x) = f_1(x), \: \pa_t u(0, x) = f_2(x),
\end{equation}
where $0 \leq q(t, x) \in C^{\infty}(\R_t \times \R^3),\:q(t, x) = 0$ for $|x| \geq \rho > 0$
and
\begin{equation} \label{eq:1.2}
\sup_{t \in \R, |x|\leq \rho}|\pa_t^k\pa_{x}^{\alpha} q(t,x)|\leq C_{k, \alpha},\: \forall k, \forall \alpha.
\end{equation}
Set 
$$\|u(t, x)\|_{\hc} = \|u(t, x)\|_{H^1(\r3)} + \|u_t(t, x)\|_{L^2(\r3)}.$$ 
For the Cauchy problem for the linear operator $\pa_t^2u - \Delta_x u+ q(t, x) u $ there exist potentials $q(t, x) \geq 0$ periodic in time with period $T > 0$ such that for suitable initial data $f = (f_1, f_2) \in \hc(\r3) = H^1(\r3) \times L^2(\r3)$ we have
$$\|u(t, x) \|_{H^1(\r3)} \geq C e^{\alpha |t|}$$
with $C > 0, \: \alpha > 0$ (see \cite{CPR}, \cite{PT}). This phenomenon is related to the so called parametric resonance.
On the other hand, adding a nonlinear term $u^3$ for the Cauchy problem (\ref{eq:1.1}) there are no parametric resonances and for every potential $q$  the solution $u(t, x)$ is defined globally for $t \in \R$ and satisfies a polynomial bound
$$\|u(t, x)\|_{H^1(\r3)} \leq B_1 (1 + B_0 |t|)^2$$
with constants $B_0 > 0, B_1 > 0$ depending on $q$ and the initial data $f \in \hc$. This result has been obtained in Theorem 2, \cite{PT} and the proof was based on the inequality
$$X'(t) \leq C X(t)^{1/2},$$
where  
$$X(t)=\frac{1}{2}\int_{\R^3}\big(|\partial_t u|^2 + |\nabla_x u|^2 + q |u|^2 +\frac{1}{2} |u|^{4}\big)dx.$$ 
If fact, the local Strichartz estimates and Theorem 2 in \cite{PT} hold for every non-negative potential $q(t, x) \in C^{\infty}(\R_t \times \R^3)$ with compact support with respect to $x$ satisfying the estimates (\ref{eq:1.2})
since in the proofs of these results the periodicity of $q$ with respect to $t$ is not used.\\
 
 In this paper we study the problem (\ref{eq:1.1}) with initial data $f \in H^k(\r3) \times H^{k-1}(\r3),\: k \geq 2.$
Throughout the paper we consider Cauchy problems with real-valued initial data $f$ and real-valued solutions. 
First in Section 2 we establish a local result and we show the existence and uniqueness of solution for $t \in [s, s+ \tau_k]$ with initial data $f \in H^k(\r3) \times H^{k-1}(\r3)$ on $t = s$ and
$$\tau_k = c_k(1 + \|(f_1, f_2)\|_{\hc(\r3)})^{-\gamma}, \: \gamma > 0,$$
where $c_k$ depends on $q$ and $k$ (see Proposition 2.1). It is important to notice that $\tau_k$ depends on the norm $\|f\|_{\hc}$ and since we have a global bound for the $\hc$ norm of $(u, u_t)(t, x)$, the interval  of local existence depends on the $\hc$ norm of the initial data.
We prove this result without using local Strichartz estimates. Next we show that the global solution in $\R$ is in $H^k(\r3)$ for 
all $t \in \R$ and the problem is to examine if the norm $\|u(t, x)\|_{H^k(\r3)}, \:k \geq 2,$ is polynomially bounded.
To do this, it is not possible to define a suitable energy $Y_k(t) \geq 0$ involving
$$\int_{\r3} (\|u(t, x)\|^2_{H^k(\r3)} + \|u_t(t, x)\|_{H^{k-1}(\r3)}^2)dx$$
for which $Y_k'(t) \leq C_k Y_k^{\gamma_k}(t),\: 0 <\gamma_k < 1.$ To overcome this difficulty, we follow another argument 
based on Lemma A.1 (see Appendix) which has an independent interest and apply local Strichartz estimates for the nonlinear equation. We study first the case $k = 2$ in Section 4 and by induction one covers the case $k \geq 3$ in Section 5. Our principal result is the following
\begin{thm} For every potential $q$ and every $k \geq 2$ the problem $(\ref{eq:1.1})$ with initial data $f \in H^k(\r3) \times H^{k-1}(\r3)$ has a global solution $u(t, x)$ and there exist $A_k > 0$ and $m_k \geq 2$ depending on $q, \:k$ and $\|f\|_{\hc}$ such that
\begin{equation} \label{eq:1.3}
\|u(t, x)\|_{H^k(\r3)} +\|\pa_t u(t, x)\|_{H^k(\r3)} \leq A_k (1 + |t|)^{m_k}, \: t \in \R.
\end{equation}
\end{thm}

We refer to \cite{PTV} and the references therein for other results about polynomial bounds for the solutions of Hamiltonian partial differential equations. The method of the proof of Theorem 1.1 basically follows the approach in \cite{PTV}. The main difficulty compared to \cite{PTV} is that in our situation, we do not have uniform bound on the $H^1(\r3)$ norm and for that purpose we need to apply the estimate of Lemma A.1 in the Appendix.
 
\section{Existence and uniqueness of local solutions in $H^k(\R^3),\: k\geq 3$}
In this section we study the existence and uniqueness of local solutions of the Cauchy problem
\begin{equation} \label{eq:2.1}
\begin{cases}   u_{tt} - \Delta_ x u + q(t, x) u + u^3 = 0, t \in [s,s+ \tau],\:x \in \R^3,\\
u(s, x) = f_1(x), u_t(s, x) = f_2(x),\end{cases}
\end{equation}
where $f = (f_1, f_2) \in H^k(\r3) \times H^{k-1}(\r3),\: k\geq 1, 0 < \tau < 1.$ We assume that $[s, s+ \tau] \subset [0, a],$ where $a > 1$ is fixed. The cases $k = 1, 2$ have been investigated in Section 3, \cite{PT} by using the norms
$$\|u\|_{S_{k-1}}: = \|(u,u_t)\|_{C([s,s+ \tau],H^k(\R^3) \times H^{k-1}(\R^3))}.$$
For $k=1$ the space $S_0$ has been denoted as $S.$ The number $\tau$ is given by 
\begin{equation} \label{eq:2.2}
\tau = c_1(1 + \|(f_1, f_2)\|_{\hc})^{-\gamma} < 1
\end{equation}
 with some positive constants $c_1 > 0, \gamma > 0$ depending on $q$. The case $k \geq 3$ can be handled by a similar argument and we will show that with $\tau =\tau_k$ defined by (\ref{eq:2.2}) with the constant $c_1$ replaced by $0 < c_k \leq c_1$ depending on $k$ and $q$, one has a local existence and uniqueness in the interval $[s, s+ \tau_k].$ Consider the linear problem
\begin{equation} \label{eq:2.3}
 \begin{cases}\partial_t^2 u_{n+1} - \Delta u_{n+1} + q(t, x) u_{n+1}+u_n^3= 0,\: n \geq 0,\\
 u_{n+1}(s, x) = f_1(x), \: \partial_t u_{n+1}(s, x) = f_2(x)\end{cases}
\end{equation}
for $t \in [s, s+ \tau_k]$ with $u_0 = 0.$ For the solution of the above problem  with right hand part $- u_n^3$ and $f = (f_1, f_2)$ we have a representation
\begin{equation} \label{eq:2.4}
(u_{n+1}, (u_{n+1})_t) = U_0(t-s)f - \int_s^t \Bigl[U_0(t - \tau)Q(\tau) u_{n+1}(\tau, x) + U_0(t - \tau)Q_0 u_n^3(\tau, x)\Bigr] d\tau.
\end{equation}
Here  $ U_0(t,s): \hc \rightarrow \hc$ is the propagator related to the free wave equation in $\R^3$ (see Section 2, \cite{PT}) and 
$$Q(\tau) = \begin{pmatrix} 0 & 0\\      
q(\tau, x) & 0 \end{pmatrix},\: Q_0= \begin{pmatrix} 0 & 0\\ 1 & 0\end{pmatrix}.$$
To estimate $\|u_{n+1}\|_{S_{k}},$ we apply the operator 
$$L_k =\begin{pmatrix}(1- \Delta)^{k/2} & 0\\ 0 & ( 1- \Delta)^{(k- 1)/2}\end{pmatrix}.$$
Notice that this operator commute with $U_0(t- \tau)$ and $\|U_0(t-s)\|_{\hc \to \hc}  \leq A$
for $|t- s| \leq 1$ with $A > 0$ independent on $k$. Therefore
$$\|U_0(t- s)L_k f\|_{\hc} \leq C \|f\|_{H^{k+1} \times H^k}$$
and
$$\Bigl\|\int_s^{\tau} U_0(t- \tau) L_k Q(\tau)u_{n+1}(\tau, x) d\tau\Bigr \|_{\hc} \leq \int_s^{\tau}\| U_0(t-\tau)L_k Q(\tau) u_{n+1}\|_{\hc}d\tau \leq A_k \tau \|u_{n+1}\|_{S_k}.$$
 For $ A_k \tau\leq 1/2 $ with $A_k >0$, depending on $k$ and $q$, the term involving $Q(\tau)u_{n+1}$ in (\ref{eq:2.4})
can be absorbed by $\|u_{n+1}\|_{S_k}$ and we deduce
$$\|u_{n+1}\|_{S_k} \leq C\|(f_1, f_2)\|_{H^{k+1}(\R^3) \times H^k(\R^3)} + C \|u_n^3\|_{L^1([s, s+ \tau], H^{k}(\R^3))}.$$
Here and below the constants $C$ depend on $k$ and $q$ and they may change from line to line but we will omit this in the notations.
Next we define the norm 
$$\|f\|_{H^{s, p}(\r3)}: = \|(1- \Delta_x)^{s/2}f\|_{L^p(\R^3)},\: 1 < p \leq\infty.$$
We will use the following {\it product estimate}
\begin{equation} \label{eq:2.5}
\|fg\|_{H^{s, p}} \leq A_{s,p} \|f\|_{L^{q_1}} \|g\|_{H^{s,q_2}} + A_{s,p} \|g\|_{L^{r_1}} \|f\|_{H^{s, r_2}},
\end{equation}
provided 
$$\frac{1}{p} = \frac{1}{q_1} + \frac{1}{q_2} = \frac{1}{r_1}+\frac{1}{r_2},\: q_1, r_1\in (1, \infty],\: q_2, r_2 \in (1, \infty].$$
For the proof of the classical estimate (\ref{eq:2.5}) we refer to \cite{T}. We apply (\ref{eq:2.5}) with $p = 2, q_1 = 3, q_2= 6, r_1 =6, r_2=3$ and get
$$\|u_n^3\|_{H^{k}(\r3)} \leq C \|u_n\|_{H^{k, 6}(\r3)}\|u_n\|_{L^6(\r3)}^2 + C \|u_n^2\|_{H^{k, 3}(\r3)} \|u_n\|_{L^6(\r3)}.$$
For the term involving $u_n^2$ in the above inequality we apply the same estimate with $p = 3, q_1= q_2= r_1= r_2 = 6$ and deduce
$$\|u_n^2\|_{H^{k, 3}(\r3)} \leq 2C \|u_n\|_{H^{k, 6}(\r3)} \|u_n\|_{L^6(\r3).}$$
Consequently, by Sobolev embedding theorem
$$\|u_n^3\|_{H^{k}(\r3)} \leq C_1 \|u_n\|_{H^{k+1}(\r3)} \|\nabla_x u_n\|_{L^2(\r3)}^2.$$
This implies
$$\int_s^{s + \tau} \|u_n^3(t, x)\|_{H^{k}(\r3)} dt \leq C_1\tau \|u_n\|_{L^{\infty}([s, s+ \tau]), H^1(\R^3))}^2\|u_n\|_{S_{k}}.$$
On the other hand, for the solution $u_n$ we have the estimate 
$$\|u_n\|_{C([s,s+\tau], H^1(\R^3))} \leq 2C_0 \|(f_1, f_2)\|_{\hc}, \: \forall n \geq 1$$
with some constant $C_0 >0$ depending on $q$ (see Section 3, \cite{PT}) and we deduce the bound
$$ C\|u_n^3\|_{L^1([s, s+ \tau], H^k(\R^3))} \leq C C_1\tau (2C_0)^2 \|(f_1, f_2)\|_{\hc}^2\|u_n\|_{S_{k+1}}.$$
Thus choosing 
$2CC_1\tau (2C_0)^2 \|(f_1, f_2)\|^2_{\hc} \leq 1,$
we may prove by induction the estimate 
\begin{equation} \label{eq:2.6}
\|u_n\|_{S_{k}}\leq 2C \|(f_1,f_2)\|_{H^{k+1}(\R^3) \times H^k(\R^3)},\: \forall n \geq 1.
\end{equation}
Next, let $w_n =u_{n+1} - u_n$ be a  solution of the problem
 $$\partial_t^2 w_n - \Delta w_n + q(t, x) w_n =  u_{n-1}^3 - u_{n}^3, \,\, w_n(0, x) = \partial_t w_n(0, x) = 0.$$
By using the inequality
$$\Bigl|v^3- w^3  \Bigr|\leq  2 |v-w|\Bigl(|v|^2 + |w|^2\Bigr),$$
 we can similarly show that 
\begin{equation*}
\|u_{n+1}-u_n\|_{S_k}\leq \frac{1}{2}\|u_{n}-u_{n-1}\|_{S_k}
\end{equation*}
which implies the convergence of $(u_n)_{n \geq 0}$ with respect to the $\|\cdot\|_{S_k}$ norm.
 Repeating the argument of Section 3, \cite{PT}, we obtain local existence and uniqueness. Thus we get the following
\begin{prop} For every $k \geq 1$ there exist $C_k > 0,\:c_k > 0$ and $\gamma > 0$ depending on $q$ and $k$ such that for every $(f_1,f_2) \in H^k(\r3) \times H^{k-1}(\r3)$ there is a unique solution $(u,u_t) \in C([s, s + \tau_k], H^k(\r3) \times H^{k-1}(\r3)$ of the problem $(\ref{eq:2.1})$ on $[s, s + \tau_k]$ with $\tau_k = c_k( 1 + \|(f_1, f_2)\|_{\hc})^{-\gamma}$. Moreover, the solution satisfies
\begin{equation} \label{eq:2.7}
\|u\|_{S_{k}} \leq C_k \|(f_1, f_2)\|_{H^{k}(\r3) \times H^{k_1}(\r3)}.
\end{equation}
\end{prop}

 It is important to note that 
for every $k,$ $\tau_k$ depends on the $\hc$ norm of the initial data.

In \cite{PT} it was proved that one has a global solution $(u, u_t) \in C(\R, \hc(\r3))$ with initial data $(f_1, f_2) \in \hc(\r3).$ It is natural to expect that for $(f_1, f_2) \in H^{k}(\r3) \times H^{k-1}(\r3)$ we have a global solution $(u, u_t) \in C(\R,H^k(\r3) \times H^{k-1}(\r3)).$ 

Let $a >1$ be fixed and let $k \geq 1.$ We wish to prove that the global solution with initial data $f \in H^{k+1}(\r3) \times H^{k}(\r3)$ is such that
\begin{equation} \label{eq:2.8}
(u, u_t) (t, x) \in H^{k +1}(\r3)\times H^{k}(\r3),\: 0 \leq t \leq a.
\end{equation}
According to Theorem 2 in \cite{PT}, for $0 \leq t \leq a$ we have an estimate
$$\|(u, u_t)(t, x)\|_{\hc} \leq B_a = \|f\|_{\hc}+ a(B_1 + B_2 a),$$
where $B_1 > 0$ and $B_2 > 0$ depend only on $\|f\|_{\hc}.$ Consider
\begin{equation} \label{eq:2.9}
\tau_k(a) = c_k(1 + B_a)^{-\gamma}.
\end{equation}
First for $0\leq t \leq \tau_k(a)$ we apply Proposition 2.1. Next we apply Proposition 2.1 for the problem with initial data on $t = \frac{2}{3}\tau_k(a)$ which is bounded by (2.7). Thus we obtain a solution in $[0, \frac{5}{3}\tau_k(a)]$ and we continue this procedure by step $\frac{2}{3} \tau_k(a)$. On every step the norm $H^{k+1}(\r3) \times H^{k}(\r3)$ of $(u, u_t)$ will increase with a constant $C_k$.
Finally, if 
$$\frac{3}{2} a \leq m \tau_k(a) \leq \frac{3}{2} (a +1),$$
 we deduce
\begin{eqnarray} \label{eq:2.10}
\|(u, u_t) (a, x)\|_{ H^{k +1} \times H^{k}} \leq C_k^m \|(f_1, f_2)\|_{H^{k+1} \times H^{k}} \nonumber \\
\leq e^{\frac{3 (a+1)}{2\tau_k(a)}\log C_k}\|(f_1, f_2)\|_{H^{k +1} \times H^{k}}.
\end{eqnarray}
Hence, we established (\ref{eq:2.8}) and one has a bound of $H^{k+1}\times H^k$ norm. Since $a$ is arbitrary, we obtain (\ref{eq:2.8}) for $t \in \R$ and a global existence for $t \in \R.$  In Section 5 we will improve (\ref{eq:2.10}) to  polynomial bounds of the Sobolev norms.

\section{Local Strichartz estimate for the nonlinear wave equation}

Our purpose is to establish a local Strichartz estimate for the solution of the Cauchy problem
\begin{equation} \label{eq:3.1}
\begin{cases}  u_{tt} - \Delta_ x u + q(t, x) u + u^3 = 0, t \in ]s, s + \tau], x \in \R^3,\\
u(s, x) = f_1(x), u_t(s, x) = f_2(x),\end{cases}
\end{equation}
where $f = (f_1, f_2) \in H^2(\r3) \times H^1(\r3),\: 0 < \tau \leq 1.$ It well known (see Proposition 1, \cite{PT}) that for the solution of the Cauchy problem
\begin{equation} \label{eq:3.2}
\begin{cases}  v_{tt} - \Delta_ x v + q(t, x) v= F, \: (t , x) \in ]s, s + \tau] \times \R^3,\\
v(s, x) = h_1(x), v_t(s, x) = h_2(x),\end{cases}
\end{equation}
we have an estimate
$$\|v(t, x)\|_{L^p([s, s + \tau], L^{r}_x(\r3))} \leq C\Bigl( \|(h_1, h_2)\|_{H^1(\r3) \times L^2(\r3)} + \|F\|_{L^1([s, s + \tau],L^2(\r3))} \Bigr),$$
where $\frac{1}{p} + \frac{3}{r} = \frac{1}{2},\: 2 < p \leq \infty.$
We will choice later $r = \frac{4 + 2 \ep}{\ep}$ with $0 < \ep \ll 1$  and this determines the choice of $p > 2.$
For the solution of (\ref{eq:3.1}) we get
\begin{eqnarray} \label{eq:3.3}
\|u(t, x) \|_ {L^p([s, s + \tau], L^{r}_x(\r3))} \leq C(p, r)\Bigl( \|u(s, x), u_t(s, x)\|_{H^1(\r3) \times L^2( \r3)} \nonumber\\
 + \tau \|u (t, x)\|^3_{L^{\infty}([s, s + \tau], H^1(\r3))}\Bigr),
\end{eqnarray}
where we have used the estimate
$$\|u^3(t, x)\|_{L^1([s, s + \tau],L^2(\r3))}   \leq  \tau \|u (t, x)\|^3_{L^{\infty}([s, s + \tau], H^1(\r3))}.$$
Next, for the solution $u(t, x) \in H^1(\r3)$ of (\ref{eq:3.1}) in $]0, s+ \tau]$ with initial data $f = (u, u_t)(0, x) \in \hc(\r3)$ we have a polynomial bound (see Section 3, \cite{PT})
$$\sup_{t \in [0, s+ \tau]}\|u(t, x)\|_{H^1(\r3)} \leq \|f\|_{\hc(\r3)} + s(B_1 +B_2 s)^2,$$
where $B_1 > 0, B_2 >0$ depend only on $\|f\|_{\hc},$
and this implies
\begin{equation} \label{eq:3.4}
\|u(t, x) \|_ {L^p([s, s + \tau], L^{r}_x(\r3))}  \leq C_1(p, r, \|f\|_{\hc})(1 + s)^6.
\end{equation}

Now we will examine the continuous dependence on the initial data of the local solution to (\ref{eq:2.1}) given in Section 2. Let $g_n = ((g_n)_1, (g_n)_2) \in H^{k+1}(\r3) \times H^k(\r3)$ be a sequence converging in $H^k(\r3) \times H^{k-1}(\r3)$ to $f = (f_1, f_2) \in H^k(\r3) \times H^{k-1}(\r3).$
Let 
$$w_n(t, x) \in C([s, s+ \tau], H^{k+1}(\r3)) \cap C^1([s, s+ \tau], H^{k}(\r3))$$
be the local solution of (\ref{eq:3.1}) with initial data $g_n$. Setting $v_n = w_n - u$, we obtain for $v_n$ the equation
$$\pa_t^2 v_n - \Delta_x v_n + q(t, x) v_n = u^3 - w^3_n.$$
By the local Strichartz estimates for the linear equation with respect to $v_n$, we get
\begin{eqnarray} \label{eq:3.5}
\|(v_n, (v_n)_t)\|_{C([s, s+ \tau], H^{k}(\r3) \times H^{k-1}(\r3))} + \|v_n\|_{L_t^{\infty}([s, s+\tau],H^{k-1, 6}_x(\r3))}\nonumber\\
\leq C_k(a) \|g_n - f\|_{H^{k}(\r3) \times H^{k-1}(\r3)} + C_k(a) \|u^3 - w_n^3\|_{L^1_t([s, s+ \tau], H^{k-1}_x(\r3))}.
\end{eqnarray}
This estimate for $k = 1, 2$ has been proved in Proposition 1, \cite{PT}. The proof for $k \geq 3$ follows the same argument. The constant $C_k(a) > 0$ depends on $k$ and on the interval $[0, a]$, where $[s, s+ \tau]\subset [0,a].$ We will omit in the notations below the dependence of the constants on $k$ and $a$.
Applying (\ref{eq:2.5}), we have
$$\|u^3 - w_n^3\|_{H^{k-1}} \leq C\|v_n\|_{H^{k-1,6}} \|u^2 + u w_n + w_n^2\|_{L^3}+ C\|v_n\|_{L^6}\|u^2 + u w_n + w_n^2\|_{H^{k-1,3}}$$
$$\leq 2 C\|v_n\|_{H^{k-1, 6}}\Bigl(\|u\|_{L^6}^2  +\|w_n\|_{L^6}^2\Bigr)
+ C\|v_n\|_{L^6}\Bigl(2\|u\|_{H^{k-1, 6}}\|u\|_{L^6}$$
$$+ 2\|w_n\|_{H^{k-1,6}}\|w_n\|_{L^6}
 + \|u\|_{H^{k-1, 6}} \|w_n\|_{L^6} + \|w_n\|_{H^{k-1,6}}\|u\|_{L^6}\Bigr) = P_n + Q_n.$$
To handle $P_n$, notice that $L^{\infty}([s, s+ \tau], L^6(\r3))$ norms of $u$ and $w_n$ by local Strichartz estimates can be estimated by $\|f\|_{\hc}$ and $\|g_n\|_{\hc}$. Therefore, for $n \geq n_0$ we have
$$\Bigl|\int_s^{s + \tau} P_n dt \Bigr| \leq A_k \tau \|v_n\|_{L^{\infty}([s, s+\tau], H^{k-1, 6}(\r3))} $$
with a constant $A_k$ depending on $C_k(a)$ and $\|f\|_{\hc}$. Hence, we may absorb $P_n$ by the left hand side of (\ref{eq:3.5}) choosing $0 < \tau \leq \frac{1}{2A_k}$ small. The analysis of $Q_n$ is easy since we proved in subsection 3.2, \cite{PT} that for all $t \in [s, s+ \tau]$ we have $\|\nabla_x v_n(t, x)\|_{L^2(\r3)} \to 0$ as $n \to \infty$ and the term in the braked $\Bigl(...\Bigr)$ for $t \in [0, a]$ is uniformly bounded with respect to $n$ according to the analysis in Section 2
and estimate (\ref{eq:2.10}).
Finally, we conclude that 
\begin{equation} \label{eq:3.6}
\|(v_n, (v_n)_t)\|_{C([s, s+ \tau], H^{k}(\r3) \times H^{k-1}(\r3))} \rightarrow_{n \to \infty} 0.
\end{equation}

\section{Polynomial bound of the $H^2(\r3)$ norm of the solution}

Let $(u(t, x), u_t(t,x)) \in C([s, s + \tau], H^2(\r3)) \times C([s, s+ \tau], H^1(\r3)),$
where $u(t, x)$ is the solution for  $t \in [s, s + \tau]$ of the Cauchy problem (\ref{eq:2.1}).
Taking the derivative $\partial_{x_j} = \partial_j, \: j = 1, 2,3,$ and noting $u_ j = \pa_j u,\: u_{j t} = \pa_{j}\pa_t u,$ one gets
in the sense of distributions
\begin{equation} \label{eq:4.1}
(u_{j t})_{t} - \Delta_ x u_j + (\pa_j q) u + q u_j + 3 u^2 u_ j = 0.
\end{equation}
It is easy to see that
$$ (\pa_j q) u + q u_j + 3 u^2 u_ j  \in C([s, s + \tau], \l2).$$
In fact, our assumption implies that $u(t, x) \in C([s, s + \tau], L_x^{\infty}(\r3))$ and this yields $u^2 u_j \in C([s, s + \tau], \l2).$ Therefore
$$(u_{j t})_{t} - \Delta_ x u_j  \in C([s, s + \tau], \l2).$$
 Multiplying the equality (\ref{eq:4.1}) by $u_{j t}$, we have
\begin{eqnarray} \label{eq:4.2}
\int \Bigl((u_{j t})_{t} - \Delta_ x u_j \Bigr) u_{j t} dx = - \int (\pa_j q) u u_{j t} dx - \int q u_j u_{j t} dx \nonumber \\
 - 3 \int u^2 u_j u_{ j t} dx  = I_1(t) + I_ 2(t) + I_ 3(t).
\end{eqnarray}
Assuming $(u(t, x), u_t(t, x))  \in C([s, s + \tau], H^3(\r3) \times H^2(\r3))$, we can write
$$I_2(t) = -\frac{1}{2} \int q \pa_t( u_j^2) dx = - \frac{1}{2} \pa_t \Bigl( \int q u_j^2 dx\Bigr) + \frac{1}{2} \int q_t u_j^2 dx,$$
$$I_3(t) = -\frac{3}{2} \int u^2 \pa_t( u_j^2) dx = - \frac{3}{2} \pa_t \Bigl(\int  u^2 u_j^2 dx \Bigr) + 3 \int u u_t u_j^2 dx.$$

After an integration by parts in the integral
$\int \Delta_x (u_j) u_{j t} dx$ 
for solutions $(u(t, x), u_t(t, x))  \in C([s, s + \tau], H^3(\r3) \times H^2(\r3))$ the equality (\ref{eq:4.2}) can be written as
\begin{eqnarray} \label{eq:4.3}
 \frac{1}{2}\pa_t \sum_{j = 1}^3 \Bigl[\int \Bigl( (u_{j t})^2 + |\nabla_x (u_j)|^2  + 3  u^2 u_j ^2  + q u_j ^2 \Bigr)(t, x) dx\Bigr] =  -\sum_{j= 1}^3\int (\pa_j q) u u_{j t} dx \nonumber \\
 + 3 \sum_{j= 1}^3\int u u_t u_j^2 dx + \frac{1}{2}\sum_{j = 1} ^3 \int q_t u_j ^2 dx= I_1(t) + J_1 (t)+ J_2(t),
\end{eqnarray}
where the derivative with respect to $t$ of the left hand side is taken in sense of distributions. 
\subsection{Justification of (\ref{eq:4.3}) for $(u(t, x), u_t(t, x)) \in C([s,s+ \tau], H^2(\R^3) \times H^1(\R^3))$}

Introduce
$$X(t): = \frac{1}{2} \sum_{j = 1}^3\int \Bigl( (u_{j t})^2 + |\nabla_x (u_j)|^2  + 3  u^2 u_j ^2  + q u_j ^2 \Bigr)(t, x) dx.$$
Notice that the function $X(t)$ is well defined. For the integral of $u^2 u_j^2$ we have

\begin{equation} \label{eq:4.4}
\int u^2 u_j^2 dx  \leq\|u\|^2_{L^4(\r3)} \|u_j\|_{L^4(\r3)}^2\leq \|u\|^{1/2}_{L^2}\|\nabla_x u\|_{L^2}^{3/2}\|u_j\|_{L^2}^{1/2} \|\nabla_x u_j\|_{L^2}^{3/2}.
\end{equation}

Also a similar argument shows that the right hand side of (\ref{eq:4.3}) is well defined and it is a continuous function of $t.$ For example,
\begin{equation} \label{eq:4.5}
\Bigl|\int u u_t u_j^2(t, x) dx \Bigr|\leq \|u_j(t,x)\|_{L^6(\R^3)}^2 \|u(t, x)\|_{L^6(\R^3)}\|u_t(t,x)\|_{L^2(\R^3)}.
\end{equation}
This implies that the derivative with respect to $t$ is taken in classical sense. Now let $(g_n, h_n) \in H^3(\r3) \times H^2(\r3)$ converge to $(u(s, x), u_t(s, x))$ in $H^2(\r3) \times H^1(\r3)$ as $n \to \infty.$ Denote as in Section 3 by $w_n(t, x)$ the local solution of (\ref{eq:3.1}) with initial data $(g_n, h_n).$ Therefore for $t \in [s, s+ \tau]$ we have
$$\int w_n^2 ((w_n)_j)^2(t, x) dx \rightarrow_{n \to \infty} \int u^2u_j^2(t, x) dx,$$
$$\int w_n (w_n)_t ((w_n)_j)^2(t, x) dx \rightarrow_{n \to \infty} \int u u_t u_j^2(t, x) dx.$$
To justify these limits, we apply the estimates (\ref{eq:4.4}) and (\ref{eq:4.5}). For example,
$$\Bigl|\int w_n (w_n)_t ((w_n)_j)^2(t, x) dx\Bigr| \leq \Bigl|\int (w_n - u) (w_n)_t ((w_n)_j)^2 dx \Bigr| $$
$$+ \Bigl| \int u ((w_n)_t - u_t) ((w_n)_j)^2 dx\Bigr| + \Bigl|\int u u_t ((w_n)_j)^2 - u_j^2) dx\Bigr|$$ 
and we use (\ref{eq:3.6}) for $k= 2$. Passing in limit in the equality (\ref{eq:4.3}) for $w_n$, we obtain it for $u$.\\

Consequently,  after an integration with respect to $t$ in (\ref{eq:4.3}), one deduces
$$X(s + \tau) = X(s) + 2\int_{s}^{s+ \tau} \Bigl(J_1(t) + J_2(t) + I_1(t)\Bigr) dt.$$

 \subsection{Estimation of $\int_{s}^{s +\tau}J_1(t)dt$}
Let $0 < \epsilon \ll1$ be a small number. First by the generalized H\"older inequality one estimates
$$|J_1(t)| \leq  3\sum_{j=1}^3\|u(t, x)\|_{L^{r}(\R^3)} \|u_t(t, x)\|_{L^{2 + \epsilon}(\R^3)}  \|u_j (t, x)\|^2_{L^4(\R^3)}$$
$$\leq 3\sum_{j=1}^3\|u(t, x)\|_{L^{r}(\R^3)} \|u_t(t, x)\|_{L^{2 + \ep}(\R^3)} \|u_j(t, x)\|^{1/2}_{L^2(\R^3)}\|u_j(t, x)\|^{3/2}_{L^6(\R^3)},$$
where
$$\frac{1}{r} = \frac{\ep}{4 + 2 \ep}.$$
According to the estimate (\ref{eq:2.7}), for $s \leq t \leq s + \tau$ by the local existence of a solution of (\ref{eq:3.1}) with initial data  $ (u(s, x), u_t(s, x)) \in H^2(\r3) \times H^1(\r3)$ on $t = s$, we  obtain
$$\|u_j(t, x)\|^{3/2}_{L^6(\r3)} \leq \|\nabla_x u_j(t, x)\|^{3/2}_{L^2(\r3)} \leq C_2 \Bigl ( \|u(s, x)\|_{H^2(\r3)} + \|u_t(s, x)\|_{H^1(\r3)}\Bigl)^{3/2}$$
with constant $C_2 > 0$ depending on $q$.
Next
$$\|u(s, x)\|^2 _{H^2(\r3)} \leq C \Bigl(\sum_{i, j = 1}^3 \|\pa_{x_i} \pa_{x_j}u(s, x)\|^2_{L^2(\r3)} + \|u(s, x)\|^2_{H^1(\r3)}\Bigr),$$
$$\|u_t(s, x)\|^2_{H^1(\r3)} \leq C\Bigl( \sum_{j= 1}^3\|u_{j t}(s, x)\|^2_{\l2} + \|u_t(s, x)\|^2_{\l2}\Bigr).$$
Notice that we have a polynomial bound with respect to $s$ for the norms $\|u(s, x)\|_{H^1(\r3)}$ and $\|u_t(s, x)\|_{L^2(\r3)}$ of the solution $u(s, x)$ (see Theorem 2, \cite{PT}).  Consequently, we obtain
 $$\sup_{t \in [s, s + \tau]} \|u_j(t, x)\|^{3/2}_{L^6(\r3)} \leq C_1\Bigl( X(s)^{3/4} +(1 + s)^3\Bigr),\:\sup_{t \in [s, s + \tau]}\|u_j(t, x)\|_{L^2(\r3)} \leq C_0( 1+ s),$$
where $C_0 > 0, C_1 > 0$ depend on $\|u(0, x)\|_{H^1(\r3)}.$\\

Now we pass to the estimate of $\|u_t(t, x)\|_{L^{2+\ep}(\r3)}.$ By H\"older inequality we obtain
$$ \Bigl| \int u_t^{2 +\ep} dx  \Bigr| =  \Bigl| \int u_t^{2(1 - \frac{\ep}{4})} u_t^{\frac{3\ep}{2}} dx  \Bigr| \leq \|u_t\|_{L^2(\r3)}^{2(1 - \ep/4)} \|u_t\|_{L^6(\r3)}^{\frac{3\ep}{2}}$$
$$ \leq  C_3(1 + t)^2 \|\nabla_x u_t\|_{L^2(\r3)}^{\frac{3\ep}{2}} \leq C_4(1 + s)^2\Bigl( X(s)^{\frac{3\ep}{4}} + (1 +s)^{3 \ep}\Bigr).$$
Hence, one deduces
$$\sup_{t \in [s, s+\tau]} \Bigl| \int u_t ^{2 +  \ep} dx \Bigr| ^{\frac{1}{2 + \ep}} \leq C_5(1+ s)^{3/2}\Bigl(X(s)^{\frac{3\ep}{8 + 4\ep}} + 1\Bigr).$$

Taking into account the above estimates, for the integral with respect to $t$ one applies the H\"older inequality and for small $\ep$ we have
$$\Bigl |\int_{s}^{s + \tau} J_1(t)dt \Bigr | \leq C_6 \tau^{1/p'} (1 + s ) ^6\|u(t, x)\|_{L^p([s, s+ \tau]; L_x^{r}(\r3))} \Bl X(s)^{\frac{3}{4} +\frac{3\ep}{8}} + 1\Br ,$$
where 
$$\frac{1}{p} + \frac{3\ep}{4 + 2 \ep} = \frac{1}{2},\:  \frac{1}{p'} +\frac{1}{p} = 1.$$
To complete the analysis, we apply the Strichartz estimate (\ref{eq:3.4}) and deduce
$$\|u(t, x)\|_{L^p([s, s+ \tau]; L_x^{r}(\r3))} \leq C(\ep) (1 + s)^6.$$

Finally for $0 < \tau \leq 1$ with  $y = 12$ we have
\begin{equation} \label{eq:4.6}
\Bigl |\int_{s}^{s + \tau} J_1(t)dt \Bigr| \leq C'(\ep)  \Bl X(s)^{\frac{3}{4} + \frac{3 \ep}{8}}+ 1 \Br(1 + s)^y.
\end{equation}
\subsection{Estimation of $\int_{s}^{s +\tau}I_1(t)dt$} We apply a similar argument.
$$|I_1(t)| \leq C\sum_{j= 1}^3 \|u(t, x)\|_{L^2(\r3)} \|u_{jt}(t, x)\|_{\l2} \leq C_7 (1 + |t|)^2 \sum_{j= 1}^3 \|u_{jt}(t, x)\|_{\l2}.$$
By the local existence result for $t \in [s, s + \tau]$ one has
$$\|u_{jt}(t, x)\|_{\l2} \leq C  ( \|u(s, x)\|_{H^2(\r3)} + \|u_t(s, x)\|_{H^1(\r3)})$$
and repeating the above argument, we deduce
\begin{equation} \label{eq:4.7}
\Bigl |\int_{s}^{s + \tau}  I_1(t)dt \Bigr | \leq C_8(X(s)^{1/2}  + 1)(1 + s)^2.
\end{equation}
\subsection{Estimation of $\int_{s}^{s +\tau}J_2(t)dt$}  This term is easy to be bounded  since we have a polynomial estimate 
$$\int u_j^2(t, x) dx \leq C_0( 1 + |t|)^2$$
and this yields
\begin{equation} \label{eq:4.8}
\Bigl |\int_{s}^{s + \tau} J_2(t)dt \Bigr | \leq C_9(1 + s)^2.
\end{equation}
Combining (\ref{eq:4.6}), (\ref{eq:4.7}), (\ref{eq:4.8}), finally we get
\begin{equation} \label{eq:4.9}
X(s + \tau) \leq X(s) + C_{10}\Bl  X(s)^{\frac{3}{4} + \frac{3\ep}{8}} + 1\Br (1 + s)^y.
\end{equation}
\subsection{Growth of $H^2(\r3)$ norm} Let $a > 1$ be a fixed number. According to \cite{PT} and Proposition 2.1, there exists a solution in
$[s , s + \tau(a)] \subset [0, a]$ with initial data $g \in H^2(\r3) \times H^1(\r3)$ on $t = s$. Here
$$\tau(a) = c \Bigl(( 1 + \|f\|_{H^1(\r3) \times \l2} + a (B_1 + B_2 a)\Bigr)^{-\gamma}< 1,$$
where $c > 0, \gamma > 0, B_1> 0, B_2 > 0$ are independent on $a$ and $f$.
 We choose $N(a) \in \N$ so that $a- \tau(a)  < N(a)\tau(a) \leq a.$
Setting $X(n\tau(a)) = \alpha_n,\:n \leq N(a),$ and exploiting (\ref{eq:4.9}), one deduces
$$\alpha_{n} \leq \alpha_{n-1} + C_{10}(\alpha_{n-1}^{7/8} + 1)(1 + n)^{12}. $$
We are in position to apply  Lemma A.1 in the Appendix and to obtain
$$X(N(a)\tau(a)) \leq \tilde{C} (N(a))^{104}$$
$$\leq \tilde{C} \Bl \frac{a}{c}\Br ^{104} \Bigl( 1 + \|f\|_{H^1(\r3) \times \l2} + a (B_1 + B_2 a)\Bigr)^{104\gamma}.$$
This estimate and the bound of the $H^1(\r3)$ norm of the solution $u(a, x)$ established in \cite{PT} imply a polynomial with respect to $a$ bound of  $\|u(a, x)\|_{H^2(\r3)} + \|\pa_t u(a, x)\|_{H^1(\r3)}$. This implies the statement of Theorem 1.1 for $k = 2.$
 
\section{Polynomial growth of the $H^k(\r3)$ norm of the solution.} 
To examine the growth of the $H^k(\r3)$ norm of the solution, we will proceed by induction. Assume that for $1\leq k \leq s-1, s\geq 3,$ we have
polynomial bounds
$$\|u(t, x)\|_{H^k_x(\r3)} + \|u_t(t, x)\|_{H^{k-1}_x(\r3)} \leq A_k (1 + |t|)^{m_k},\: t \in \R$$
for the global solution of the Cauchy problem of $u_{tt} - \Delta_x u +qu+ u^3 = 0$ with initial data $(f_1, f_2) \in H^k(\r3) \times H^{k-1}(\r3).$ Consider the equality 
$$\pa_t^2\pa_x^{\alpha}u - \Delta_ x (\pa_x^{\alpha}u) + \pa_x^{\alpha}(qu) + \pa_x^{\alpha} (u^3) = 0$$
 with $|\alpha| = s - 1.$ After an integration by parts which we can justify as in Section 4, we write
\begin{eqnarray} \label{eq:5.1}
\frac{1}{2}\frac{d}{dt}\int \Bigl(|\nabla_x \pa_x^{\alpha} u|^2 + |\pa_t \pa_x^{\alpha} u|^2\Bigr) dx \nonumber\\
= - \int \pa_x^{\alpha}(q u) \pa_x^{\alpha} \pa_t u dx - \int \pa_x^{\alpha} (u^3) \pa_x^{\alpha} \pa_t u dx = K_1(t) + K_2(t).
\end{eqnarray}
Clearly,
$$\Bigr|\int \Bigr(\pa_x^{\alpha} (u^3) \pa_x^{\alpha} \pa_t u\Bigr) dx\Bigr | \leq \|\pa_x^{\alpha} (u^3)\|_{L^2(\r3)}\|\pa_x^{\alpha} \pa_t u\|_{L^2(\r3)}.$$ 
Applying two times (\ref{eq:2.5}), one gets
$$\|\pa_x^{\alpha} (u^3)\|_{L^2(\r3)} \leq C\|\pa_x^{\alpha} u\|_{L^2(\r3)} \|u\|_{L^{\infty}(\r3)}^2$$
and by Sobolev theorem $\|u\|_{L^{\infty}(\r3)} \leq C \|u\|_{H^2(\r3)}.$ Thus by our assumption
$$\|\pa_x^{\alpha} (u^3(t, x))\|_{L^2(\r3)} \leq C A_{k-1} A_2^2 (1 + |t|)^{m_{k-1} + 2m_2}.$$
Therefore, using the notation of subsection 4.5 for $n \tau(a) \leq t \leq (n+1) \tau(a),$ one deduces
$$\|\pa_x^{\alpha} (u^3(t, x))\|_{L^2(\r3)} \leq C A_{k-1} A_2^2(1 + n)^{m_{k-1} + 2m_2}.$$
On the other hand, applying (\ref{eq:2.7}), one obtains
$$\|\pa_x^{\alpha} \pa_t u(t, x)\|_{L^2(\r3)} \leq C_{k}\Bigl((\|u(n\tau(a),x)\|_{H^{s}(\r3)} + \|u_t(n\tau(a),x)\|_{H^{s-1}(\r3)}\Bigr).$$
The analysis of $K_1(t)$ is easy and 
$$|K_1(t)| \leq C \|u(t, x)\|_{H^{s-1}(\r3)} \|\pa_x^{\alpha}\pa_t u(t, x)\|_{L^2(\r3)}$$
$$\leq C_k A_{k-1}(1 + n)^{m_{k-1}}(\|u(n\tau(a),x)\|_{H^{s}(\r3)} + \|u_t(n\tau(a),x)\|_{H^{s-1}(\r3)}).$$

Now define
$Y_k(t): = \|u(t, x)\|_{H^k(\r3)}^2 + \|\pa_t u(t, x)\|_{H^{k-1}(\r3)}^2$
and integrate the equality (\ref{eq:5.1}) from $n\tau_k(a)$ to $(n+1)\tau_k(a)$ with respect to $t$, where $0 < \tau_k(a) < 1$ is defined by (\ref{eq:2.9}).   Taking into account the above estimates, we have
$$Y_k((n+1)\tau_k(a)) \leq Y_k(n \tau_k(a)) + C_q A_{k-1}(1 + n)^{m_{k-1}} $$
$$+ C A_{k-1} A_2^2 (1+ n)^{m_{k-1} + 2m_2} Y_k^{1/2}(n\tau_k(a)).$$ 
Applying Lemma A.1 and repeating the argument of subsection 4.5, we obtain a polynomial bound for $Y_k(t)$ and this completes the proof of Theorem 1.1.

\section{Appendix}
\renewcommand{\thelem}{A.\arabic{lem}}
\renewcommand{\therem}{A.\arabic{rem}}
 \renewcommand{\theequation}{A.\arabic{equation}}
\setcounter{equation}{0}

The aim in this Appendix is to prove the following
\begin{lem} Let $\{\alpha_n\}_{n=0}^{\infty}$ be a sequence of non-negative numbers such that with some constants $0 < \gamma <1$,
$C > 0$ and $y \geq 0$ we have
$$\alpha_n \leq \alpha_{n-1} + C ((\alpha_{n-1})^{1 - \gamma}+ 1)(1 +  n)^y,\: \forall n \geq 1.$$
Then there exists a constant $\tilde{C} > 0$ such that
\begin{equation} \label{eq:A.1}
\alpha_n \leq \tilde{C}(1 + n)^{\gy},\:\forall n \geq 1.
\end{equation}
\end{lem}
\begin{rem} A similar estimate has been established in \cite{PTV} for sequences $\{\alpha_n\}$ satisfying the inequality
$$\alpha_n \leq \alpha_{n-1} + C \alpha_{n-1}^{1 - \gamma}.$$
\end{rem}
\begin{proof}We can choose a large constant $C_1 > 0$ such that
$$(\alpha_{n-1})^{1 - \gamma} + 1 \leq C_1 (\alpha_{n-1} + 1)^{1-\gamma}, \: \forall n\geq 1.$$
This implies with a new constant $C_2 > 0$ the inequality
$$\alpha_n  + 1\leq \alpha_{n-1} + 1+ C_2(\alpha_{n-1} + 1)^{1 - \gamma} (1 +  n)^y,\: \forall n \geq 1.$$
Setting $\beta_n = \alpha_n + 1,$ we reduce the proof to a sequence $\alpha_n$ satisfying the inequality 
$$\alpha_n \leq \alpha_{n-1} + C_2 (\alpha_{n-1})^{1 - \gamma}(1 +  n)^y, \:n \geq 1.$$

We will prove (\ref{eq:A.1}) by recurrence. Assume that (\ref{eq:A.1}) holds for $n-1$. Therefore
$$\alpha_n \leq \tilde{C} n^{\gy} + C_2 \Bigl( \tilde{C} n^{\gy}\Bigr)^{1-\gamma} (1+n)^y$$
$$= \tilde{C} n^{\gy} \Bigl[ 1 + C_2 \tilde{C}^{-\gamma}n^{-1 - y}(1+ n)^y\Bigr]$$
$$ = \tilde{C} (1+n)^{\gy}\Bigl( 1 -\frac{1}{n + 1}\Bigr)^{\gy}\Bigl[ 1 + C_2 \tilde{C}^{-\gamma}n^{-1}\Bigl(\frac{n}{n+1}\Bigr)^{- y}\Bigr].$$
To establish (\ref{eq:A.1}) for $n$, it is sufficient to show that for large $\tilde{C}$ one has
\begin{equation} \label{eq:A.2}
f(n): = \Bigl( 1 -\frac{1}{n + 1}\Bigr)^{\gy}\Bigl[ 1 + C_2 \tilde{C}^{-\gamma}n^{-1}\Bigl(\frac{n}{n+1}\Bigr)^{- y}\Bigr] \leq 1,\: n\geq 1.
\end{equation}
Setting $C_2 \tilde{C}^{-\gamma} = \epsilon$, a simple calculus yields
$$f'(n) = \frac{1 + y}{\gamma} \Bigl(1 - \frac{1}{n+1}\Bigr)^{\gy - 1}\frac{1}{(n+1)^2}\Bigl[ 1 + \frac{\epsilon}{n}\Bigl(\frac{n}{n+ 1}\Bigr)^{- y}\Bigr]$$
$$ + \epsilon\Bigl(1 - \frac{1}{n+1}\Bigr)^{\gy}\Bigl[ -\frac{1}{n^2}\Bigl(\frac{n}{n+ 1}\Bigr)^{-y}-y n^{-1} \frac{1}{(n+1)^2}  \Bigl(1 - \frac{1}{n+1}\Bigr)^{-y -1}\Bigr]$$
$$= \Bigl(1 - \frac{1}{n+1}\Bigr)^{\gy -1}\frac{1}{(n+ 1)^2}\Bigl[ \gy + \frac{\epsilon}{n}\gy \nn^{-\gamma} - [\epsilon \frac{n+1}{n}+ \frac{\epsilon y}{n}] \nn^{-y}\Bigr].$$
Notice that since $\frac{1}{2} \leq 1- \frac{1}{n+ 1}$, we have
$$\nn^{-\gamma} \leq \Bigl(\frac{1}{2}\Bigr)^{-\gamma}$$
which implies
$$\gy- \epsilon[\frac{n+1 + y}{n}] \nn^{-y} \geq \gy- \epsilon[\frac{n+1 + y}{n}] \Bigl(\frac{1}{2}\Bigr)^{-y}.$$
For small $\epsilon > 0$ the right hand side of the above inequality is positive. Consequently, for the derivative we have
$f'(n) > 0$ and one deduces
$$f(n) < \lim_{n \to +\infty} f(n)= 1$$
This completes the proof of (A.2).
\end{proof}

\section*{Acknowledgments} We would like to thank the referee for his/her useful comments.

\end{document}